\documentclass[12pt,draft]{article}

\usepackage[cp1251]{inputenc}
\usepackage[T2A]{fontenc}
\usepackage[english,ukrainian]{babel}
\usepackage{amsmath,amsfonts,amssymb}
\usepackage{geometry}

\begin{document}

\begin{center}

\large

\textbf{ELLIPTIC PROBLEMS\\ WITH ROUGH BOUNDARY DATA\\ IN NIKOLSKII SPACES}

\medskip

\textbf{A.A. Murach, I.S. Chepurukhina}

\medskip

\normalsize

Institute of Mathematics of the National Academy of Sciences of Ukraine, Kyiv

\medskip

E-mail: murach@imath.kiev.ua, Chepurukhina@gmail.com

\bigskip\bigskip

\large

\textbf{ЕЛІПТИЧНІ ЗАДАЧІ\\ З ГРУБИМИ КРАЙОВИМИ ДАНИМИ\\ У ПРОСТОРАХ НІКОЛЬСЬКОГО}

\medskip

\textbf{О.О. Мурач, І.С. Чепурухіна}

\normalsize

\end{center}

\medskip

\noindent We investigate a general elliptic problem given in a bounded Euclidean domain with boundary data in Nikolskii spaces of low, specifically, negative order. The right-hand side of the elliptic differential equation is supposed to be an integrable function. We establish the Fredholm property of the problem, the maximal regularity and a priori estimate of its generalized solutions in the spaces indicated. We give an application of these results to some elliptic problems with boundary conditions induced by a Gaussian white noise.

\medskip

\noindent\textbf{Keywords}: elliptic boundary-value problem, Nikolskii space, Fredholm operator, regularity of solution, a priori estimate, white noise.


\bigskip\bigskip

\noindent \textbf{Вступ.} У теорії багатовимірних крайових задач окрему увагу привертають задачі з грубими крайовими даними, тобто такими, що не є інтегровними функціями на межі (див., наприклад, [\ref{LionsMagenes71}, \ref{Roitberg96}, \ref{MikhailetsMurach14}, \ref{BehrndtHassiSnoo20}]). До останніх належать функції зі степеневими особливостями, дельта-функції та їх похідні, розподіли, породжені різними випадковими процесами. Вони  інтерпретуються у сенсі теорії узагальнених функцій (розподілів) як елементи деяких функціональних просторів від'ємного порядку. Дослідження таких задач є істотно більш складним порівняно з випадком достатньо  регулярних крайових даних, які належать до просторів додатного порядку. Тут зразу постає питання про інтерпретацію слідів узагальнених функцій низької регулярності на межі області, де розглядається задача. Як правило, такі сліди означають за допомогою граничних переходів, апроксимуючи узагальнені функції в області досить гладкими функціями у спеціально підібраних нормах. Уведення таких норм накладає деякі умови регулярності на праву частину диференціального рівняння [\ref{LionsMagenes71}], або приводить до того, що розв'язок задачі не можна інтерпретувати як узагальнену функцію в області [\ref{Roitberg96}]. Додаткові складнощі виникають, коли використовують простори розв'язків, у яких множини гладких функцій не є щільними, як, наприклад, простори Нікольського або простори Гельдера--Зігмунда низького порядку.

Мета цієї роботи~--- встановити теореми про характер розв'язності (нетеровість), максимальну регулярність і апріорну оцінку розв'язків еліптичних задач з крайовими даними у просторах Нікольського низького (у тому числі від'ємного) порядку. Вона мотивована важливими застосуваннями останніх у теорії білого шуму. Саме у термінах просторів Нікольського від'ємного порядку вдалося точно (за порядком простору) охарактеризувати регулярність гауссового білого шуму на торі [\ref{Veraar11}]. Це дозволяє нам отримати точний результат про регулярність розв'язків деяких еліптичних задач з крайовими умовами, породженими таким шумом.

\textbf{1. Постановка задачі.} Нехай $\Omega$~--- обмежена область у евклідовому просторі $\mathbf{R}^{n}$, де $n\geq2$. Припустимо, що її межа $\Gamma$ є нескінченно гладким замкненим многовидом вимірності $n-1$ (причому $C^{\infty}$-структура на $\Gamma$ породжена $\mathbf{R}^{n}$).
Розглянемо в $\Omega$ еліптичну крайову задачу вигляду
\begin{gather}\label{12f1}
Au=f\quad\mbox{в}\;\,\Omega,\\
B_{j}u=g_{j}\quad\mbox{на}\quad\Gamma,\quad j=1,\ldots,l. \label{12f2}
\end{gather}
Тут $A:=A(x,D)$~--- лінійний диференціальний оператор (л.д.о.) на $\overline{\Omega}:=\Omega\cup\Gamma$ парного порядку $2l\geq2$, а кожне $B_{j}:=B_{j}(x,D)$~--- крайовий л.д.о. на $\Gamma$ порядку $m_{j}\leq 2l-1$. Припускаємо, що усі коефіцієнти цих л.д.о. є нескінченно гладкими комплекснозначними функціями, заданими на $\overline{\Omega}$ і $\Gamma$ відповідно. Нагадаємо, що еліптичність крайової задачі \eqref{12f1}, \eqref{12f2} означає, що л.д.о. $A$ є правильно еліптичним на $\overline{\Omega}$, а набір $B:=(B_1,\ldots,B_l)$ крайових л.д.о. задовольняє умову Лопатинського щодо $A$ на $\Gamma$.

Нехай дійсне число $p>1$. Досліджуємо властивості розв'язку $u$ крайової задачі \eqref{12f1}, \eqref{12f2} у ситуації, коли крайові дані $g_{1},\ldots,g_{l}$ належать відповідним просторам Нікольського $B^{\sigma}_{p,\infty}(\Gamma)$ низької регулярності, зокрема, коли їх порядок $\sigma<0$. При цьому припускаємо, що права частина $f$ еліптичного рівняння \eqref{12f1} належить до $L_{p}(\Omega)$, тобто функція $|f|^{p}$ інтегровна на $\Omega$ (відносно міри Лебега). Слідуючи [\ref{Triebel92}, п.~2.3.1] і [\ref{Triebel86}, п.~3.2.2], нагадаємо означення  просторів Нікольського довільного дійсного порядку $\sigma$. У~статті розглядаються комплексні лінійні функціональні простори.

Довільно виберемо функцію $\varphi_0\in C^{\infty}(\mathbf{R}^{n})$ таку, що $\varphi_0(y)=1$ за умови $|y|\leq1$ та $\varphi_0(y)=0$ за умови $|y|\geq2$, де $y\in\mathbf{R}^{n}$. Для кожного номера $k\geq1$ означимо функцію  $\varphi_k(y):=\varphi_{0}(2^{-k}y)-\varphi_{0}(2^{1-k}y)$ аргументу $y\in\mathbf{R}^{n}$. Функції $\varphi_{k}$, де ціле $k\geq0$, утворюють нескінченно гладке розбиття одиниці на евклідовому просторі $\mathbf{R}^{n}$.

Нехай $\sigma$~--- довільне дійсне число. За означенням, лінійний простір $B^{\sigma}_{p,\infty}(\mathbf{R}^{n})$ складається з усіх повільно зростаючих розподілів $w$ на $\mathbf{R}^{n}$ таких, що
$$
\|w,B^{\sigma}_{p,\infty}(\mathbf{R}^{n})\|:=
\sup_{0\leq k\in\mathbf{Z}}2^{k\sigma}
\biggl(\;\int\limits_{\mathbf{R}^{n}}
|F^{-1}[\varphi_{k}Fw]|^{p}(x)dx\biggr)^{1/p}<\infty,
$$
де $F$ і $F^{-1}$ позначають відповідно пряме і обернене перетворення Фур'є. Цей простір є повним відносно норми  $\|w,B^{\sigma}_{p,\infty}(\mathbf{R}^{n})\|$. Він не залежить з точністю до еквівалентності норм від вказаного вибору функції $\varphi_0$ [\ref{Triebel92}, теорема 2.3.2(i)]. Простір $B^{\sigma}_{p,\infty}(\mathbf{R}^{n})$ був уведений (в інший еквівалентний спосіб) і досліджений С.М.~Нікольським для  $\sigma>0$ (див. монографії [\ref{Nikolskii77}, с.~151] і [\ref{Triebel80}, с.~201] та наведені там посилання). Зауважимо, що С.М.~Нікольський вибрав позначення $H^{\sigma}_{p}(\mathbf{R}^{n})$ для цього простору (яке часто застосовують для просторів Соболєва). Позначення, яке використовуємо ми, загальноприйняте і свідчить про те, що шкала просторів Нікольського є специфічною частиною класу просторів Бєсова $B^{\sigma}_{p,q}(\mathbf{R}^{n})$.

Простори Нікольського на $\Omega$ і $\Gamma$ будуються за простором $B^{\sigma}_{p,\infty}(\mathbf{R}^{n})$ у стандартний спосіб (див., наприклад, [\ref{Triebel86}, п.~3.2.2]). А саме, лінійний простір $B^{\sigma}_{p,\infty}(\Omega)$ складається, за означенням, із звужень в область $\Omega$ усіх розподілів $w\in B^{\sigma}_{p,\infty}(\mathbf{R}^{n})$ і наділений нормою
$$
\|u,B^{\sigma}_{p,\infty}(\Omega)\|:=
\inf\bigl\{\|w,B^{\sigma}_{p,\infty}(\mathbf{R}^{n})\|:w\in B^{\sigma}_{p,\infty}(\mathbf{R}^{n}),\;\,u=w\;\,\mbox{в}\;\,\Omega\bigr\},
$$
де $u\in B^{\sigma}_{p,\infty}(\Omega)$. Простір $B^{\sigma}_{p,\infty}(\Omega)$ є повним відносно цієї норми.

Простір $B^{\sigma}_{p,\infty}(\Gamma)$ складається з усіх розподілів на $\Gamma$, які в локальних координатах на $\Gamma$ дають елементи простору $B^{\sigma}_{p,\infty}(\mathbf{R}^{n-1})$. Дамо детальне означення. З~$C^{\infty}$-структури на многовиді $\Gamma$ довільно виберемо скінченний набір локальних карт $\pi_j:\mathbf{R}^{n-1}\leftrightarrow\Gamma_{j}$, де $j=1,\ldots,\lambda$, таких, що відкриті множини $\Gamma_{1},\ldots,\Gamma_{\lambda}$ утворюють покриття $\Gamma$. Крім того, виберемо функції $\chi_j\in C^{\infty}(\Gamma)$, де $j=1,\ldots,\lambda$, які утворюють розбиття одиниці на $\Gamma$ і задовольняють умову $\mathrm{supp}\,\chi_j\subset\Gamma_j$. (Як звичайно, $\mathrm{supp}\,\chi_j$ позначає замикання в топології $\Gamma$ множини усіх точок $x\in\Gamma$ таких, що $\chi_j(x)\neq0$.) За означенням, лінійний простір $B^{\sigma}_{p,\infty}(\Gamma)$ складається з усіх розподілів $h$ на $\Gamma$ таких, що $(\chi_{j}h)\circ\pi_{j}\in B^{\sigma}_{p,\infty}(\mathbf{R}^{n-1})$ для кожного номера $j\in\{1,\ldots,\lambda\}$ і наділений нормою
$$
\|h,B^{\sigma}_{p,\infty}(\Gamma)\|:=
\sum_{j=1}^{\lambda}\,\|(\chi_{j}h)\circ\pi_{j},
B^{\sigma}_{p,\infty}(\mathbf{R}^{n-1})\|;
$$
тут $(\chi_{j}h)\circ\pi_{j}$ позначає зображення розподілу $\chi_{j}h$ у локальній карті $\pi_{j}$. Простір $B^{\sigma}_{p,\infty}(\Gamma)$ повний відносно цієї норми і з точністю до еквівалентності норм не залежить від вказаного вибору локальних карт і розбиття одиниці на $\Gamma$ [\ref{Triebel86}, твердження 3.2.3(ii)].

Зауважимо, що множини $C^{k}(\overline{\Omega})$ і $C^{k}(\Gamma)$, де ціле $k\geq0$, лежать відповідно у просторах $B^{\sigma}_{p,\infty}(\Omega)$ і $B^{\sigma}_{p,\infty}(\Gamma)$ порядку $\sigma\leq k$, але не є щільними у них. Останнє випливає з [\ref{Triebel80}, теорема 2.3.2(a)].

\textbf{2. Основні результати.} Нехай, як і раніше, $1<p<\infty$. Еліптична крайова задача \eqref{12f1}, \eqref{12f2} має таку фундаментальну властивість у просторах Нікольського додатного порядку:
відображення $u\mapsto(Au,Bu)$, де $u\in B^{s}_{p,\infty}(\Omega)$, є нетеровим обмеженим оператором на парі просторів
\begin{equation}\label{AB-positive}
(A,B):B^{s}_{p,\infty}(\Omega)\to B^{s-2l}_{p,\infty}(\Omega)\times
\prod_{j=1}^{l}B^{s-m_j-1/p}_{p,\infty}(\Gamma)
\quad\mbox{для кожного}\;\,s>2l.
\end{equation}
Крім того, ядро цього оператора лежить у $C^{\infty}(\overline{\Omega})$ і разом з індексом не залежить від $s$ і~$p$. Для регулярних еліптичних задач ця властивість доведена у монографії [\ref{Triebel80}, теорема 5.5.2(a)]. Для загальних еліптичних задач, розглянутих нами, вона міститься у результатах робіт [\ref{Murach94UMJ}, п.~3, теорема] і [\ref{Johnsen96}, теорема~5.2]. Позначимо через $N$ та $\alpha$ відповідно ядро та індекс оператора \eqref{AB-positive}.

Принагідно нагадаємо, що лінійний обмежений оператор $T:E_{1}\rightarrow E_{2}$, де $E_{1}$ і $E_{2}$~--- повні нормовані простори, називають нетеровим, якщо його ядро $\{x\in E_{1}:Tx=0\}$ і коядро $E_{2}/T(E_{1})$ скінченновимірні. Якщо цей оператор нетерів, то його область значень замкнена в $E_{2}$ і він має скінченний індекс, який, за означенням, дорівнює різниці вимірностей ядра і коядра.

Стосовно оператора \eqref{AB-positive} зауважимо, що образ $Au\in B^{s-2l}_{p,\infty}(\Omega)$ означений у сенсі теорії розподілів для довільного $u\in B^{s}_{p,\infty}(\Omega)$, де $s$~--- будь-яке дійсне число. Крім того, образ $B_{j}u\in B^{s-m_j-1/p}_{p,\infty}(\Gamma)$ означений у сенсі теореми про сліди для кожного вказаного $u$, якщо $s>m_j+1/p$ (див. [\ref{Triebel80}, теорема 4.7.1(b)]); тут фіксовано номер $j\in\{1,\ldots,l\}$. Проте, якщо $s\leq m_j+1/p$, то образ $B_{j}u$ не означений коректно для довільного $u\in B^{s}_{p,\infty}(\Omega)$ навіть як деякий розподіл на $\Gamma$. А~саме, у цьому випадку відображення $u\mapsto\nobreak B_{j}u$, де $u\in C^{\infty}(\overline{\Omega})$, не допускає продовження до неперервного лінійного оператора на парі просторів $B^{s}_{p,\infty}(\Omega)$ і $D'(\Gamma)$, що випливає з [\ref{Triebel80}, теорема 4.7.1(d)]. Як звичайно, $D'(\Gamma)$ позначає лінійний топологічний простір усіх розподілів на~$\Gamma$.

Отже, лінійний обмежений оператор \eqref{AB-positive} не означається коректно для довільного $s<2l$. У цьому випадку замість  $B^{s}_{p,\infty}(\Omega)$ беремо вужчий простір
$$
B^{s}_{p,\infty}(A,L_p,\Omega)=\bigl\{u\in B^{s}_{p,\infty}(\Omega):Au\in L_p(\Omega)\bigr\}
$$
як область визначення оператора $(A,B)$. Цей простір наділений нормою графіка $\|u,B^{s}_{p,\infty}(\Omega)\|+\|Au,L_p(\Omega)\|$ і є повним відносно неї. Для довільного $u\in B^{s}_{p,\infty}(A,L_p,\Omega)$ коректно означений розподіл $B_{j}u\in D'(\Gamma)$ у сенсі, який вказано нижче.

Розглянемо лінійний простір
$$
S'(A,L_p,\Omega):=\{u\in S'(\Omega):Au\in L_{p}(\Omega)\},
$$
де $S'(\Omega)$~--- лінійний топологічний простір звужень в область $\Omega$ усіх повільно зростаючих розподілів на $\mathbf{R}^{n}$. Оскільки множина $\Omega$ обмежена, то $S'(A,L_p,\Omega)$ є об'єднанням усіх просторів $B^{s}_{p,\infty}(A,L_p,\Omega)$, де $s\in\mathbf{R}$. При цьому $B^{s}_{p,\infty}(A,L_p,\Omega)$ неперервно (але не щільно) вкладений у $B^{\sigma}_{p,\infty}(A,L_p,\Omega)$, якщо $\sigma<s$. У лінійному просторі $S'(A,L_p,\Omega)$ означаємо збіжність за таким правилом: говоримо, що $u_{k}\to u$ в $S'(A,L_p,\Omega)$ (при $k\to\infty$), якщо
$u_{k}\to u$ в $B^{s}_{p,\infty}(A,L_p,\Omega)$ для деякого дійсного~$s$.
Множина $C^{\infty}(\overline{\Omega})$ щільна у просторі $S'(A,L_p,\Omega)$ відносно цієї збіжності, що випливає з [\ref{MikhailetsMurach14}, теорема 4.25(i)] з огляду на теорему вкладення [\ref{Triebel80}, теорема 4.6.2]. Для довільного $u\in S'(A,L_p,\Omega)$ виберемо послідовність функцій $u_{k}\in C^{\infty}(\overline{\Omega})$ таку, що $u_{k}\to u$ в $S'(A,L_p,\Omega)$. Згідно з [\ref{MikhailetsMurach14}, теорема 4.25(ii)] послідовність $(B_{j}u_{k})_{k=1}^{\infty}$ має границю у топологічному просторі $D'(\Gamma)$ і ця границя не залежить від вказаного вибору послідовності функцій $u_{k}$. Останню границю і беремо як $B_{j}u$. У результаті є коректно означеним лінійний неперервний оператор
\begin{equation}\label{AB-extended}
(A,B):S'(A,L_p,\Omega)\to L_p(\Omega)\times(D'(\Gamma))^{l}.
\end{equation}

\textbf{Теорема 1.} \it Нехай $s\leq2l$ і $1<p<\infty$. Тоді звуження відображення \eqref{AB-extended} на простір $B^{s}_{p,\infty}(A,L_p,\Omega)$ є обмеженим оператором на парі просторів
\begin{equation}\label{12f4}
(A,B):B^{s}_{p,\infty}(A,L_p,\Omega)\to
L_p(\Omega)\times\prod_{j=1}^{l}B^{s-m_j-1/p}_{p,\infty}(\Gamma).
\end{equation}
Оператор \eqref{12f4} нетерів з ядром $N$ та індексом $\alpha$ (які не залежать від $s$ і $p$). \rm

Теорему 1 доречно порівняти з результатом статті [\ref{Murach94UMJ}, п.~3] про нетеровість еліптичних крайових задач у просторах Нікольського--Ройтберга довільного дійсного порядку~$s$. Ці простори служать областю визначення обмеженого оператора, породженого еліптичною крайовою задачею \eqref{12f1}, \eqref{12f2}, і нетеровість якого встановлена у вказаній статті. Вони складаються з елементів, які не є функціями чи розподілами, якщо $s<2l-1/p$. Простори крайових даних для цього оператора такі самі як і в~\eqref{12f4}, а простір правих частин складається у випадку $s<2l$ з усіх розподілів $w\in B^{s-2l}_{p,\infty}(\mathbf{R}^{n})$ таких, що $\mathrm{supp}\,w\subset\overline{\Omega}$. Отже, вказаний результат знаходиться поза межами теорії розподілів в області~$\Omega$ у випадку, який нас цікавить.

Дослідимо локальні (впритул до частини межі області $\Omega$) властивості узагальнених розв'язків еліптичної крайової задачі \eqref{12f1}, \eqref{12f2}. Розподіл $u\in S'(\Omega)$ називаємо узагальненим розв'язком цієї задачі з правими частинами $f\in L_p(\Omega)$ і $g_{1},\ldots,g_{l}\in D'(\Gamma)$, якщо $Au=f$ у сенсі теорії розподілів і $B_{j}u=g_{j}$ для кожного номера $j\in\{1,\ldots,l\}$ у вказаному вище сенсі для $u\in S'(A,L_p,\Omega)$. З умови $Au\in L_p(\Omega)$ та еліптичності л.д.о. $A$ випливає, що $\chi u\in B^{s}_{p,\infty}(\Omega)$ для кожного числа $s<2l$ і довільної функції $\chi\in C^{\infty}(\overline{\Omega})$, яка дорівнює нулю в околі межі $\Gamma$. Втім, якщо  $\chi(x)\neq0$ на деякій частині межі, то включення $\chi u\in B^{s}_{p,\infty}(\Omega)$ буде виконуватися лише за певних умов на крайові дані $g_{1},\ldots,g_{l}$ в околі множини $\Gamma\cap\mathrm{supp}\,\chi$. Якщо $s\geq2l$, то, звісно, потрібна ще деяка умова на $f$ в околі $\mathrm{supp}\,\chi$, щоб гарантувати це включення. Уведемо простори, у термінах яких зручно формулювати ці умови.

Нехай $U$ --- відкрита підмножина простору $\mathbf{R}^{n}$ така, що $\Omega_{0}:=\Omega\cap U\neq\varnothing$ і $\Gamma_{0}:=\Gamma\cap U\neq\varnothing$. Позначимо через $B^{\sigma,\mathrm{loc}}_{p,\infty}(\Omega_{0},\Gamma_{0})$, де $\sigma\in\mathbf{R}$, лінійний простір усіх розподілів $u\in\mathcal{S}'(\Omega)$ таких, що $\chi u\in B^{\sigma}_{p,\infty}(\Omega)$ для кожної функції $\chi\in C^{\infty}(\overline{\Omega})$, яка задовольняє умову $\mathrm{supp}\,\chi\subset\Omega_0\cup\Gamma_{0}$.
Аналогічно, позначимо через $B^{\sigma,\mathrm{loc}}_{p,\infty}(\Gamma_{0})$
лінійний простір усіх розподілів $h\in D'(\Gamma)$ таких, що $\chi h\in B^{\sigma}_{p,\infty}(\Gamma)$ для кожної функції $\chi\in C^{\infty}(\Gamma)$, яка задовольняє умову   $\mathrm{supp}\,\chi\subset\Gamma_{0}$.

\textbf{Теорема 2.} \it Нехай $s\in\mathbf{R}$ і $1<p<\infty$. Припустимо, що розподіл $u\in S'(\Omega)$ є узагальненим розв'язком еліптичної крайової задачі \eqref{12f1}, \eqref{12f2}, праві частини якої задовольняють умови
\begin{equation}\label{f-cond}
f\in L_p(\Omega)\cap B^{s-2l,\mathrm{loc}}_{p,\infty}(\Omega_{0},\Gamma_{0})
\end{equation}
і $g_j\in B^{s-m_j-1/p,\mathrm{loc}}_{p,\infty}(\Gamma_0)$ для кожного $j\in\{1,\ldots,l\}$. Тоді $u\in B^{s,\mathrm{loc}}_{p,\infty}(\Omega_{0},\Gamma_{0})$. \rm

\textbf{Зауваження 1.} Якщо $s\leq2l$, то умова \eqref{f-cond} еквівалентна включенню $f\in L_p(\Omega)$ згідно з [\ref{Triebel80}, теорема 4.6.1(a,b)]. Отже, у цьому випадку, висновок теореми~2 стає таким: $u\in B^{s,\mathrm{loc}}_{p,\infty}(\Omega,\Gamma_{0})$.

Цю теорему доповнює локальна апріорна оцінка узагальненого розв'язку~$u$.

\textbf{Теорема 3.} \it Нехай $s\in\mathbf{R}$ і $1<p<\infty$. Припустимо, що розподіл $u\in S'(\Omega)$ задовольняє умови теореми~$2$. Довільно виберемо число $r>0$ і функції $\chi,\eta\in C^{\infty}(\overline{\Omega})$ такі, що $\mathrm{supp}\,\chi\subset\mathrm{supp}\,\eta\subset
\Omega_{0}\cup\Gamma_{0}$ і $\eta=1$ в околі $\mathrm{supp}\,\chi$. Тоді
\begin{equation}\label{12f6}
\begin{gathered}
\|\chi u,B^{s}_{p,\infty}(\Omega)\|\leq
c\,\biggl(\|\eta f,L_{p}(\Omega)\|+
\|\eta f,B^{s-2l}_{p,\infty}(\Omega)\|+\\
+\sum_{j=1}^{l}\|\eta g_{j},B^{s-m_j-1/p}_{p,\infty}(\Gamma)\|+
\|\eta u,B^{s-r}_{p,\infty}(\Omega)\|\biggr),
\end{gathered}
\end{equation}
де $c$~--- деяке додатне число, яке не залежить від $u$, $f$ і $g_{1},\ldots,g_{l}$. \rm

\textbf{Зауваження 2.} Якщо $s\leq2l$, то в апріорній оцінці \eqref{12f6} можна прибрати доданок $\|\eta f,B^{s-2l}_{p,\infty}(\Omega)\|$, а якщо $s>2l$, то~--- доданок $\|\eta f,L_{p}(\Omega)\|$. Це випливає з [\ref{Triebel80}, теорема 4.6.1(a,b)].

\textbf{Зауваження 3.} Зворотною до теорем 2 і 3 є така локальна властивість крайових даних у просторах Нікольського. Нехай $s\in\mathbf{R}$ і $1<p<\infty$. Припустимо, що розподіл $u\in B^{s,\mathrm{loc}}_{p,\infty}(\Omega_{0},\Gamma_{0})$ є узагальненим розв'язком крайової задачі \eqref{12f1}, \eqref{12f2}, де
$f\in L_p(\Omega)$ і $g_{1},\ldots,g_{l}\in D'(\Gamma)$. Тоді
$g_j\in B^{s-m_j-1/p,\mathrm{loc}}_{p,\infty}(\Gamma_0)$ для кожного $j\in\{1,\ldots,l\}$, причому
$$
\|\chi g_j,B^{s-m_j-1/p}_{p,\infty}(\Gamma)\|\leq
c_{0}\bigl(\|\eta u,B^{s}_{p,\infty}(\Omega)\|+
\|\eta f,L_p(\Omega)\|\bigr).
$$
Тут функції $\chi$ і $\eta$ такі як у теоремі~3, а $c_{0}$~--- деяке додатне число, яке не залежить від $u$, $f$ і $g_{1},\ldots,g_{l}$. Ця властивість виконується і у випадку, коли набір $B$ крайових л.д.о. не задовольняє умову Лопатинського щодо $A$ на $\Gamma$.

Як бачимо, теореми 2 і 3 дозволяють за регулярністю крайових даних еліптичної задачі \eqref{12f1}, \eqref{12f2} у просторах Нікольського (як завгодно низького порядку) встановити максимальну регулярність її узагальненого розв'язку (у термінах ізотропних просторів). Теореми 1--3 варто порівняти з новітньою роботою [\ref{Hummel21JEE}], де для деяких еліптичних з параметром крайових задач у півпросторі доведено твердження [\ref{Hummel21JEE}, теорема 6.3] про існування, єдиність, регулярність і апріорну оцінку розв'язків в анізотропних просторах, побудованих на основі ряду просторів низької регулярності, заданих на межі (зокрема, просторів Нікольського). У цій роботі не досягнуто максимальної регулярності розв'язків, як зазначає її автор.

Обговоримо коротко методику доведення теорем 1--3. Теорема~1 виводиться з її аналогу для $L_{p}$-просторів Соболєва в $\Omega$ за допомогою дійсної інтерполяції $(E_{0},E_{1})_{\theta,\infty}$ нормованих просторів $E_{0}$ і $E_{1}$. При цьому використовуються інтерполяційні формули, доведені, наприклад, в [\ref{Triebel86}, п.~3.3.6], та їх версія для для простору $B^{s}_{p,\infty}(A,L_p,\Omega)$, яка обґрунтовується подібно до [\ref{KasirenkoMikhailetsMurach19}, теорема~2]. Зазначений аналог теореми~1 доведено в [\ref{LionsMagenes63VI}, теорема~8.1] для регулярних еліптичних крайових задач і просторів Соболєва цілого порядку $s<0$. Для довільних еліптичних задач, розглянутих нами, він обґрунтовується подібно до [\ref{MikhailetsMurach14}, теореми 4.25 і 4.27]. Теореми 2 і 3 доводяться у спосіб в ідейному плані близький до того, що використано у роботі [\ref{AnopDenkMurach21CPAA}, доведення теорем 4.7 і 4.13] стосовно еліптичних  задач з крайовими даними низької регулярності у гільбертових узагальнених просторах Соболєва. При цьому у доведенні теореми~3 треба скористатися згаданим вище результатом статті [\ref{Murach94UMJ}, п.~3] про нетеровість еліптичних крайових задач у просторах Нікольського--Ройтберга.

\textbf{3. Застосування.} Розглянемо застосування теореми~1 до еліптичних задач з білим шумом у крайових умовах. Нехай $(\Theta,K,\mathbf{P})$~--- ймовірнісний простір, тобто $K$~--- деяка $\sigma$-алгебра підмножин довільної множини $\Theta\neq\varnothing$, а $\mathbf{P}$~--- деяка міра на $K$, підпорядкована умові $\mathbf{P}(\Theta)=1$. За означенням (див., наприклад, [\ref{Veraar11}, с.~390]), гауссів білий шум на $\Gamma$~--- це випадкова величина $\xi:\Theta\to D'(\Gamma)$ така, що: 1) для кожного $v\in C^{\infty}(\Gamma)$ комплексна числова випадкова величина $\xi(v):\Theta\to\mathbf{C}$ має нормальний розподіл; 2) виконується рівність
$
\mathbf{M}[\xi(v_1)\overline{\xi(v_2)}]=
C\int_{\Gamma}v_1(x)\overline{v_2(x)}dS
$
для довільних функцій $v_1,v_2\in C^{\infty}(\Gamma)$ і деякого числа $C>0$, незалежного від $v_1$ і $v_2$. Тут, як звичайно, $\mathbf{M}$~--- математичне сподівання, $\xi(v)$~--- значення розподілу (узагальненої функції) $\xi$ на основній функції $v\in C^{\infty}(\Gamma)$, а $dS$~--- елемент площі поверхні~$\Gamma$.

Для гауссового білого шуму $\xi$ на торі $\mathbf{T}^{d}$ вимірності $d\geq1$ доведено в [\ref{Veraar11}, теорема~3.4], що
\begin{equation}\label{white-noise-in-B}
\mathbf{P}\{\xi\in B^{-d/2}_{p,\infty}(\mathbf{T}^{d})\}=1
\quad\mbox{для кожного}\;\;p\in(1,\infty).
\end{equation}
Тут не можна замінити простір Нікольського на більш вузькі простори Бєсова $B^{-d/2}_{p,q}(\mathbf{T}^d)$, де $p,q\in(1,\infty)$ або $p=q=\infty$, бо для них $\mathbf{P}\{\xi\in B^{-d/2}_{p,q}(\mathbf{T}^d)\}=0$. Звідси випливає, що $\mathbf{P}\{\xi\in B^{\sigma}_{p,\infty}(\mathbf{T}^d)\}=0$ як тільки $\sigma>-d/2$. Отже, у формулі \eqref{white-noise-in-B} значення параметрів $\sigma=-d/2$ і $q=\infty$ є точними, як і умова $p<\infty$.

Припустимо, що обмежена область $\Omega\subset\mathbf{R}^{d+1}$ має межу $\Gamma=\mathbf{T}^{d}$, де $d\geq1$. Розглянемо у цій області крайову задачу Діріхле для рівняння Пуассона
\begin{equation}\label{Poisson-Dirichlet-noise}
\Delta u=f\quad\mbox{в}\;\,\Omega,\qquad
\gamma_{0}u=\xi\quad\mbox{на}\;\,\mathbf{T}^{d},
\end{equation}
де $f\in L_{p}(\Omega)$ для деякого $p\in(1,\infty)$, а $\xi$~--- гауссів білий шум на торі $\mathbf{T}^{d}$. Тут $\Delta$~--- оператор Лапласа, а $\gamma_0$ позначає крайовий оператор (порядку нуль), який кожній функції $u\in C^{\infty}(\overline{\Omega})$ ставить у відповідність її слід (тобто звуження) на $\mathbf{T}^{d}$ і продовжується єдиним чином за допомогою граничного переходу до неперервного лінійного оператора на парі просторів $S'(A,L_p,\Omega)$ і $D'(\mathbf{T}^{d})$. Розглянута задача є еліптичною крайовою задачею вигляду \eqref{12f1}, \eqref{12f2}, де $l=1$, $m_1=0$, $N=\{0\}$ і $\alpha=0$.

З теореми 1 і властивості \eqref{white-noise-in-B} випливає такий результат:

\textbf{Теорема 4.} \it Для $\mathbf{P}$-майже всіх $\omega\in\Theta$ крайова задача~\eqref{Poisson-Dirichlet-noise} має єдиний узагальнений розв'язок $u(\omega,\cdot)$ класу $B^{1/p-d/2}_{p,\infty}(\Omega)$. Він задовольняє апріорну оцінку
\begin{equation*}
\|u(\omega,\cdot),B^{1/p-d/2}_{p,\infty}(\Omega)\|\leq c\,\bigl(\|f,L_{p}(\Omega)\|+
\|\xi(\omega),B^{-d/2}_{p,\infty}(\mathbf{T}^{d})\|\bigr),
\end{equation*}
де $c$~--- деяке додатне число, яке не залежить від $f$, $\xi$ і $\omega$. \rm

Цей результат є точним за порядком простору Нікольського в $\Omega$. Використання просторів Бєсова $B^{\sigma}_{p,q}$ з показником $q<\infty$ не дозволяє отримати подібний точний результат, тобто досягти граничного значення $\sigma=1/p-d/2$. У роботі [\ref{AnopDenkMurach21CPAA}, висновок~8.4] його було досягнуто для деяких гільбертових узагальнених просторів Соболєва (випадок $p=2$), але висновок теореми~4 є більш сильним з огляду на включення у них просторів Нікольського, використаних у цій теоремі.

\bigskip\bigskip

\noindent ЦИТОВАНА ЛІТЕРАТУРА

\begin{enumerate}

\item\label{LionsMagenes71}
Лионс Ж.-Л., Мадженес Э. Неоднородные граничные задачи и их приложения. Москва: Мир, 1971.

\item\label{Roitberg96}
Roitberg Ya. Elliptic boundary value problems in the spaces of
distributions. Dordrecht: Kluwer Acad. Publ., 1996.

\item\label{MikhailetsMurach14}
Mikhailets V.A., Murach A.A. H\"ormander spaces, interpolation, and elliptic problems. Berlin, Boston: De Gruyter, 2014.

\item\label{BehrndtHassiSnoo20}
Behrndt J., Hassi S., de Snoo H. Boundary value problems, Weyl functions, and differential operators. Cham: Springer, 2020.

\item\label{Veraar11}
Veraar M. Regularity of Gaussian white noise on the d-dimensional torus. \textit{Banach Center Publ.} 2011. \textbf{95}. P. 385--398.

\item\label{Triebel92}
Triebel H. Theory of function spaces. II. Basel: Birkh\"aser, 1992.

\item\label{Triebel86}
Трибель Х. Теория функциональных пространств. Москва: Мир, 1986.

\item\label{Nikolskii77}
Никольский С.М. Приближение функций многих переменных и теоремы вложения. Москва: Наука, 1977.

\item\label{Triebel80}
Трибель Х. Теория интерполяции, функциональные пространства, дифференциальные операторы. Москва: Мир, 1980.

\item\label{Murach94UMJ}
Мурач А.А. Эллиптические краевые задачи в полных шкалах пространств типа Никольского. \textit{Укр. мат. журн.} 1994. \textbf{46}, №~12. С. 1647--1654.

\item\label{Johnsen96}
Johnsen J. Elliptic boundary problems and the Boutet de Monvel calculus in Besov and Triebel-Lizorkin spaces. \textit{Math. Scand.} 1996. \textbf{79}, №~1. P. 25--85.

\item\label{Hummel21JEE}
Hummel F. Boundary value problems of elliptic and parabolic type
with boundary data of negative regularity. \textit{J. Evol. Equ.} 2021. https://doi.org/10.1007/s00028-020-00664-0

\item\label{KasirenkoMikhailetsMurach19}
Kasirenko T., Mikhailets V., Murach A. Sobolev-like Hilbert spaces induced by elliptic operators. \textit{Complex Anal. Oper. Theory.} 2019. \textbf{13}, №~3. P. 1431--1440.

\item\label{LionsMagenes63VI}
Lions J.-L., Magenes E. Probl\'emes aux limites non homog\'enes.~VI. \textit{J.~d'Analyse Math.} 1963. \textbf{11}. P. 165--188.

\item\label{AnopDenkMurach21CPAA}
Anop A., Denk R., Murach A. Elliptic problems with rough boundary data in generalized Sobolev spaces. \textit{Commun. Pure Appl. Anal.} 2021. \textbf{20}, №~2. P. 697--735.

\end{enumerate}

\bigskip

\noindent REFERENCES

\begin{enumerate}

\item
Lions, J.-L. \& Magenes, E. (1972). Non-homogeneous boundary-value problems and applications, vol.~I. Berlin: Springer.

\item
Roitberg, Ya. (1996). Elliptic boundary value problems in the spaces of distributions. Dordrecht: Kluwer Acad. Publ.

\item
Mikhailets, V.A. \& Murach, A.A. (2014). H\"ormander spaces, interpolation, and elliptic problems. Berlin, Boston: De Gruyter.

\item
Behrndt, J., Hassi S. \& de Snoo, H. (2020). Boundary value problems, Weyl functions, and differential operators. Cham: Springer.

\item
Veraar, M. (2011). Regularity of Gaussian white noise on the d-dimensional torus. Banach Center Publ., 95, pp. 385-398.

\item
Triebel, H. (1992). Theory of function spaces. II. Basel: Birkh\"aser.

\item
Triebel, H. (1983). Theory of function spaces. Basel: Birkh\"auser.

\item
Nikol'skii, S.M. (1977). Approximation of functions of several variables and imbedding theorems [2-nd edn]. Moscow: Nauka. (Russian.)

\item
Triebel, H. (1995). Interpolation theory, function spaces, differential operators [2-nd edn]. Heidelberg: Johann Ambrosius Barth.

\item
Murach, A.A. (1994). Elliptic boundary value problems in complete scales of Nikol'skii-type spaces. Ukrainian Math J., 46, No.~12, pp. 1827-1835.

\item
Johnsen, J. (1996). Elliptic boundary problems and the Boutet de Monvel calculus in Besov and Triebel-Lizorkin spaces. Math. Scand., 79, No.~1, pp. 25-85.

\item
Hummel, F. (2021). Boundary value problems of elliptic and parabolic type with boundary data of negative regularity. J. Evol. Equ., https://doi.org/10.1007/s00028-020-00664-0

\item
Kasirenko, T., Mikhailets, V. \& Murach, A. (2019). Sobolev-like Hilbert spaces induced by elliptic operators. Complex Anal. Oper. Theory, 13, No.~3, pp. 1431-1440.

\item
Lions, J.-L. \& Magenes, E. (1963) Probl\'emes aux limites non homog\'enes.~VI. J.~d'Analyse Math., 11, pp. 165-188.

\item
Anop, A., Denk, R. \& Murach, A. (2021) Elliptic problems with rough boundary data in generalized Sobolev spaces. Commun. Pure Appl. Anal.,  20, No.~2, pp. 697-735.

\end{enumerate}

\end{document}